\documentclass[11pt]{article}

\usepackage{amssymb, color}
\usepackage{amscd}

\newtheorem{theorem}{Theorem}
\newtheorem{proposition}[theorem]{Proposition}

\newtheorem{remark}[theorem]{Remark}

\def\blacksquare{\hbox{\vrule width 4pt height 4pt depth 0pt}}
\def\cqd{\ \ \ \hbox{}\nolinebreak\hfill $\blacksquare \  \  \
\  $ \par{}\medskip}

\begin{document}

\title{Factorization theorems for dominated polynomials}
\author{Geraldo Botelho\thanks{Supported by CNPq Project 202162/2006-0.}\,\,,
Daniel Pellegrino\thanks{Supported by CNPq Grant 308084/2006-3.}~ and Pilar
Rueda\thanks{Supported by MEC and FEDER Project MTM2005-08210.\hfill \newline 2000 Mathematics
Subject Classification. Primary 46G25; Secondary 47B10.}}

\date{}
\maketitle \vspace*{-1.0em}

\begin{abstract}
In this note we prove that the factorization theorem for dominated
polynomials given in \cite{BoPeRu} is equivalent to an alternative
factorization scheme that uses classical linear techniques and a
linearization process.
 However, this alternative scheme is shown not to be satisfactory until the equivalence is proved.
\end{abstract}

\section{Introduction and notation}
Let $n$ be a positive integer and $1 \leq p < +\infty$. Generalizing
the well-known Pietsch Factorization Theorem for absolutely
$p$-summing linear operators, in \cite{BoPeRu} we prove the
existence of a canonical prototype of an $n$-homogeneous
$p$-dominated polynomial such that any other $n$-homogeneous
$p$-dominated polynomial between Banach spaces is the composition of
a restriction of it with some bounded linear operator. Using the
linearization of homogeneous polynomials on projective symmetric
tensor products, in this note we obtain an alternative factorization
scheme, which is nonnatural at first glance (in a sense that will
become clear in Remark \ref{remark}). However, by relating
 the
$p$-dominated norm of a $p$-dominated polynomial with the sup
norms of the mappings involved in the factorization diagrams, we prove that this alternative approach is indeed equivalent to the factorization scheme 
provided in \cite{BoPeRu}. Although proofs are based mainly on the results in \cite{BoPeRu} and their proofs therein, the interest of this note is to clarify that the alternative factorization diagram, a priori non satisfactory, turns out to be equivalent to the previous one and that the construction in \cite{BoPeRu} comes into help to achieve this end.

From now on, $E$ and $F$ will denote (real or complex) Banach spaces, and $B_ {E^*}$ the closed
unit ball of the dual space $E^*$ of $E$. When endowed with the weak* topology, $B_ {E^*}$ is a
compact space, and, as usual, $C(B_{E^*})$ means the space of all continuous functions on
 $B_ {E^*}$ endowed with the supremum norm $\| \phantom f\|_\infty$. Let $n$ be a positive integer
 and $1\leq p< +\infty$. Given a regular Borel probability
 measure $\mu$ on $B_{E^*}$, let $j_{p/n} \colon C(B_{E^*})\longrightarrow L_{p/n}(\mu)$ be the canonical map.

Let us denote by ${\mathcal P}(^nE;F)$ the space of all continuous
$n$-homogeneous polynomials from $E$ into $F$ endowed with the
usual sup norm. For the theory of homogeneous polynomials we refer
to the excellent monograph \cite{din}.

Consider the $n$-homogeneous polynomial $j_{p/n}^n \colon
C(B_{E^*})\longrightarrow L_{p/n}(\mu)$ given by
$j_{p/n}^n(f)=j_{p/n}(f^n)$. By $e$ we mean the evaluation map
from $E$ into $C(B_{E^*})$ given by $e(x)=e_x$, where $e_x(x^*)=
\langle x^*,x \rangle$ for any $x^*\in B_{E^*}$. The restriction
of $j_{p/n}^n$ to $e(E)$ will be denoted $(j_{p/n}^e)^n$, so
$(j_{p/n}^e)^n: e(E) \longrightarrow E^{p/n}$ where
$E^{p/n}:=(j_{p/n}^e)^n\circ e(E)$.

A continuous $n$-homogeneous polynomial $P \colon E\longrightarrow
F$ is said to be $p$-dominated if there is a constant $C>0$ such
that
$$
\left( \sum_{j=1}^m\|P(x_j)\|^{p/n}\right)^{n/p}\leq C
\sup_{x^*\in B_{E^*}}\left(\sum_{j=1}^m |\langle
x^*,x_j\rangle|^p\right)^{n/p}
$$
for any positive integer $m$ and any finite sequence of vectors
$x_1, \ldots,x_m \in E$. The infimum of the constants $C$ for which
the above inequality holds is denoted by $\|P\|_{d,p}$. It is well
known that $\|\phantom P\|_{d,p}$ makes the set ${\mathcal
P}_{d,p}(^nE;F)$ of all $p$-dominated $n$-homogeneous polynomials
from $E$ into $F$ a Banach space if $p\geq n$ and a complete
$p/n$-normed space if $p< n$ (see e.g. \cite{note}). If $n=1$ we
recover the classical ideal of absolutely $p$-summing linear
operators.

\section{Alternative factorization scheme}

In order to describe  the factorization theorem for dominated polynomials proved in \cite{BoPeRu} and get the alternative approach
we need to introduce some notation. By $\otimes^{n,s}_{\pi_s}E$ we denote the $n$-fold symmetric
tensor product of $E$ endowed with the projective s-tensor norm $\pi_s$ (see \cite{klaus} for
definitions and main properties). As usual, $\hat\otimes^{n,s}_{\pi_s}E$ denotes the completed
space. Let $\delta^E_n \colon E\longrightarrow \otimes^{n,s}_{\pi_s}E$ be the diagonalization map
given by $\delta^E_n(x)=x\otimes \cdots\otimes x$. Given $P \in {\cal P}(^n E;F)$, $P^L \colon
\otimes^{n,s}_{\pi_s}E \longrightarrow F$ is the linearization of $P$, that is the unique
continuous linear map from $\otimes^{n,s}_{\pi_s}E$ to $F$ fulfilling $P^L\circ \delta^E_n=P$. It
is well known that $\|P^L\|=\|P\|$.

According to \cite{BoPeRu}, there is an injective linear operator
$\delta \colon \otimes^{n,s}_{\pi_s}E \longrightarrow C(B_{E^*})$
such that $\delta (x\otimes \cdots\otimes x)(x^*)= \langle x^*,x
\rangle^n$ for any $x\in E$ and $x^*\in B_{E^*}$, and
$j_{p/n}\circ\delta(\otimes^{n,s}_{\pi_s}E)=(j_{p/n}^e)^n\circ
e(E)=E^{p/n}$. So, $E^{p/n}$ is a linear subspace of
$L_{p/n}(\mu)$. Moreover, in \cite{BoPeRu} a norm $\pi_{p/n}$ is
defined on $E^{p/n}$ by:
$$
\pi_{p/n}((j_{p/n} \circ \delta)(\theta)) := \inf \left \{\sum_{i=1}^m |\lambda_i|\|(j_{p/n} \circ
\delta)(x_i \otimes \cdots \otimes x_i)\|_{L_{p/n}}\right\},
$$
where the infimum is taken over all representations of $\theta \in \otimes^{n,s}_{\pi_s}E$ of the
form $\theta = \sum_{i=1}^m \lambda_i x_i \otimes \cdots \otimes x_i$ with $m \in \mathbb{N}$,
$\lambda_i \in \mathbb{K}$ and $x_i \in E$.

Let us recall the factorization theorem for dominated polynomials that appears in \cite[Theorem 4.6]{BoPeRu}:

\begin{theorem}\label{factorization}{\rm \cite[Theorem 4.6]{BoPeRu}} Let $P \in {\mathcal P}(^nE;F)$ and $1 \leq p < +\infty$.
Then $P$ is $p$-dominated if and only if there is a regular Borel
probability measure $\mu$ on $B_{E^*}$ with the weak* topology and
a continuous linear operator $u \colon \left( E^{p/n},
\pi_{p/n}\right) \longrightarrow F$ such that the following
diagram commutes
\begin{center}
$\begin{CD}
E                      @ > P >>            F\\
@V{e}VV                  @ AAuA\\
e(E) @>{\left(j^e_{p/n}\right)^n}>>  E^{p/n}\\
@VVV @VVV\\
C(B_{E^*}) @>{\left(j_{p/n}\right)^n} >> L_{p/n}(\mu)
\end{CD}$
\end{center}
Moreover, $ \|u\| =\|P\|_{d,p}.
$
\end{theorem}

Although the equality $ \|u\| =\|P\|_{d,p}$ does not appear in the
statement of \cite[Theorem 4.6]{BoPeRu}, it can be obtained from
\cite{BoPeRu} as follows:

 From the proof of \cite[Theorem 4.6]{BoPeRu} we have that
 $ \|u\| \leq\|P\|_{d,p}$.
 In \cite[Lemma 4.5]{BoPeRu} it is proved that $(j^e_{p/n})^n$ is $p$-dominated and that
 for every $x_1,\ldots,x_m\in E$,
$$
\left(\sum_{i=1}^m \pi_{p/n}\left((j^e_{p/n})^n(e(x_i)) \right)^{p/n} \right)^{n/p} \leq
\sup_{\psi\in B_{C(B_{E^*})^*}}\left(\sum_{j=1}^m |\langle \psi,e(x_j) \rangle|^p\right)^{n/p}.
$$
It follows that $\|(j^e_{p/n})^n\|_{d,p} \leq 1$.
 Since $P = u\circ \left(j^e_{p/n}\right)^n \circ e$, using the ideal property \cite[Proposition
2.2]{BoPeRu} we get
\[\|P\|_{d,p} = \|u\circ (j^e_{p/n})^n \circ e \|_{d,p} \leq \|u\| \,
 \|(j^e_{p/n})^n\|_{d,p} \, \|e\|^n = \|u\|.\]

 \medskip

  By combining the classical linear factorization techniques with a linearization
process we get a different approach to the  factorization scheme given in Theorem \ref{factorization}. Naturally, $j_p^e$ denotes the restriction of $j_p$ to $e(E)$ onto its range.

\begin{proposition}\label{factorization2} Let $P \in {\mathcal P}(^nE;F)$ and $1 \leq p < +\infty$.
Then $P$ is $p$-dominated if and only if there is a regular Borel probability measure $\mu$ on
$B_{E^*}$ with the weak* topology, a closed subspace $G_{E\!,p}$ of $L_{p}(\mu)$ containing
$j_p(e(E))$ and a continuous linear operator $v \colon \otimes^{n,s}_{\pi_s}G_{E\!,p}
\longrightarrow F$ such that the following diagram commutes
\begin{center}
$\begin{CD}
E                      @ > P >>            F\\
@V{e}VV                  @ AAvA\\
e(E) @>{\delta^{G_{E\!,p}}_n\circ j^e_{p}}>>  \otimes^{n,s}_{\pi_s}G_{E\!,p}\\
\end{CD}$
\end{center}
Moreover, $\|v\|=\|P\|_{d,p}.$
\end{proposition}

\noindent {\sc Proof:}  Assume that $P$ is $p$-dominated. Adapting to polynomials the proof of the
linear case (the multilinear case can be found in \cite[Teorema 3.6]{teseDavid}) or combining
Pietsch Factorization Theorem  with \cite[Proposition 3.4]{BoPeRu}, we can conclude that there is a
regular Borel probability measure $\mu$ on $B_{E^*}$  and a continuous $n$-homogeneous polynomial
$Q$ from $G_{E\!,p}:=\overline{j_p^e\circ e (E)}\subset L_p(\mu)$ to $F$ such that the following
diagram commutes
\begin{center}
$\begin{CD}
E                      @ > P >>            F\\
@V{e}VV                  @ AAQA\\
e(E) @>{j^e_{p}}>>  G_{E\!,p}\\
@VVV @VVV\\
C(B_{E^*}) @>{j_{p}} >> L_p(\mu)
\end{CD}$
\end{center}
that is, $Q\circ j_p^e\circ e=P$. Moreover, $\|Q\|\leq \|P\|_{p,d}$. Defining $v:=Q^L$ as the
linearization of $Q$ on $\otimes_{\pi_s}^{n,s}G_{E\!,p}$ we obtain the desired commutative diagram.

As to the converse, by \cite[Proposition 3.4]{BoPeRu} we know that $\delta_n^{G_{E\!,p}}\circ
j_p^e$ is $p$-dominated, so $P=v\circ(\delta_n^{G_{E\!,p}}\circ j_p^e)\circ e$ is $p$-dominated as
well.

Denoting by $\pi_p(j_p^e)$ the $p$-summing norm of $j_p^e$, the following computation completes the
proof:
\begin{eqnarray}
\|P\|_{d,p} \!\!\!\! & = & \!\!\!\! \|Q^L \circ \delta^{G_{E\!,p}}_n \circ j_p^e \circ e\|_{d,p}
\leq \|Q^L\| \, \|  \delta^{G_{E\!,p}}_n \circ j_p^e
\|_{d,p} \, \|e\|^n\nonumber\\
\!\!\!\! & \leq & \!\!\!\!  \|Q^L\| \, \|\delta^{G_{E\!,p}}_n\| \,(\pi_p(j_p^e))^n \leq \|Q^L\| =
\|Q\| \leq \|P\|_{d,p}. \nonumber
\end{eqnarray}
 \cqd

\begin{remark}\label{remark} \rm (a) It is worth mentioning that the intermediate factorization scheme
$P=Q\circ j_p^e\circ e$ in the proof above provides the factorization of $p$-dominated polynomials
 through a canonical absolutely $p$-summing {\it linear operator.}

 \noindent (b) Observe that the $n$-homogeneous polynomial $\delta^{G_{E\!,p}}_n \circ j_p^e$ is
 $p$-dominated by \cite[Proposition 3.4]{BoPeRu}. Theorem \ref{factorization2} assures that,
 like the polynomial
$(j_{p/n}^e)^n$ of Theorem
 \ref{factorization},  $\delta^{G_{E\!,p}}_n \circ j_p^e$ is a canonical prototype of a $p$-dominated
polynomial through which every $p$-dominated polynomial factors. Let us explain why Theorem
\ref{factorization2} is a somewhat nonnatural generalization of Pietsch Factorization Theorem to
$p$-dominated polynomials. Remember that Pietsch Factorization Theorem gives the factorization of a
$p$-summing linear operator $u \colon E \longrightarrow F$ through a canonical prototype of a
$p$-summing linear operator from a subspace of $C(B_{E^*})$ into some subspace of $L_p(\mu)$, where
$\mu$ is some regular Borel probability measure on $B_{E^*}$. Observe that this is exactly what
happens for polynomials in Theorem \ref{factorization}, whereas in Theorem \ref{factorization2} we
do not know (up to now) if $\otimes^{n,s}_{\pi_s}G_{E\!,p}$ is a subspace of $L_{q_n}(\mu)$ for
some number $q_n$ with $q_1 = p$ and some regular Borel probability measure $\mu$ on
$(B_{E^*},w^*)$.
\end{remark}

We finish the note by proving that the completion $\overline E^{p/n}$ of the space $E^{p/n}$ of
Theorem \ref{factorization} is isometrically isomorphic to the completion
$\hat\otimes^{n,s}_{\pi_s}G_{E\!,p}$ of the space $\otimes^{n,s}_{\pi_s}G_{E\!,p}$ of Proposition
\ref{factorization2}. This shows that the two factorization schemes are equivalent, a fact which
reinforces the role of Theorem \ref{factorization} as a generalization of Pietsch Factorization
Theorem to dominated polynomials.

\begin{theorem}\label{tensor}
The spaces $(\overline E^{p/n},\pi_{p/n})$ and
$\hat\otimes^{n,s}_{\pi_s} G_{E\!,p}$ are isometrically
isomorphic.
\end{theorem}
\noindent {\sc Proof:} It suffices to show that $(\overline
E^{p/n},\pi_{p/n})$ and $\hat\otimes^{n,s}_{\pi_s} j_p^e\circ
e(E)$ are isometrically isomorphic. The proof is based on the proof of \cite[Proposition 4.2]{BoPeRu}.  Consider $T:=\otimes^{n,s}j_p^e\circ e$ the symmetric
$n$-fold tensor product of the linear operator $j_p^e\circ e$. By
\cite[1.7]{klaus}, $T \colon \otimes^{n,s}_{\pi_s}E
\longrightarrow \otimes^{n,s}_{\pi_s} j_p^e\circ e(E)$ is a linear
operator and satisfies
$$T(x \otimes \cdots \otimes x) = j_p^e(e(x)) \otimes \cdots \otimes j_p^e(e(x)) {\rm~for~every~} x \in E.$$
On the other hand,
$T$ is also injective and
\begin{equation}\label{iso}
\pi_s(T(\theta)) = \pi_{p/n}(j_{p/n} \circ \delta(\theta))
 \end{equation}
 for all $\theta \in \otimes^{n,s}E$.
 Then $T \colon \otimes^{n,s}_{\pi_s}E \longrightarrow
\otimes^{n,s}_{\pi_s} j_p^e\circ e(E)$ is a linear bijection.
Hence $j_{p/n}\circ \delta \circ T^{-1} \colon
\otimes^{n,s}_{\pi_s} j_p^e\circ e(E)\longrightarrow E^{p/n}$ is a
linear bijection, and by
(\ref{iso}) this map is an isometry. Taking completions we get the
desired conclusion.
 \cqd

\noindent {\sc Acknowledgements.} This note was written while G.B.
was a CNPq Postdoctoral Fellow in the Departamento de An\'alisis
Matem\'atico at Universidad de Valencia. He thanks Pilar Rueda and
the members of the department for their kind hospitality.

\noindent [Geraldo Botelho] Faculdade de Matem\'atica,
Universidade Federal de Uberl\^andia, 38.400-902 - Uberl\^andia,
Brazil,  e-mail: botelho@ufu.br.

\medskip

\noindent [Daniel Pellegrino] Departamento de Matem\'atica,
Universidade Federal da Para\'iba, 58.051-900 - Jo\~ao Pessoa,
Brazil, e-mail: dmpellegrino@gmail.com.

\medskip

\noindent [Pilar Rueda] Departamento de An\'alisis Matem\'atico,
Universidad de Valencia, 46.100 Burjasot - Valencia, Spain,
e-mail: pilar.rueda@uv.es.

%

\end{document}